\documentclass[a4paper,11pt]{article}

\usepackage{amsthm,amsmath,amssymb}
\usepackage[latin1]{inputenc}
\usepackage[english]{babel}
\usepackage{mathscinet}
\usepackage{graphics}
\usepackage{graphicx}
\usepackage[all,cmtip]{xy}
\usepackage{multirow}
\usepackage{hyperref}

\theoremstyle{plain}

\theoremstyle{definition}

\theoremstyle{remark}

\title{\bf Volumes of Solids of Revolution. \\ A Unified Approach}
\author{Jorge Mart\'{\i}n-Morales\footnote{Partially supported by the projects MTM2010-21740-C02-02, ``E15 Grupo Consolidado Geometr\'{\i}a'' from the goverment of Aragón, FQM-333 from ``Junta de Andalucía'', and PRI-AIBDE-2011-0986 Acción Integrada hispano-alemana.} \, and\, Antonio M.~Oller-Marcén \\ jorge@unizar.es, \ \ oller@unizar.es}
\date{Centro Universitario de la Defensa - IUMA.\\ Academia General Militar, Ctra.~de Huesca s/n.\\ 50090, Zaragoza, Spain.}

\begin{document}

\maketitle


\section{Introduction}
\label{intro}
The computation of the volume of solids of revolution is a very common topic in undergraduate calculus courses \cite{LYE,STE}. Usually two methods are presented in textbooks, namely:

\begin{enumerate}
\item The disk method, which roughly consists of decomposing the solid into slices that are perpendicular to the axis of revolution.
\item The shell method, that considers the solid as a series of concentric cylindrical shells wrapping the axis.
\end{enumerate}

From a geometrical point of view these two methods are quite different and it is the shape of the solid what motivates the choice between one or another. Nevertheless, from an analytical point of view, both methods are equivalent and they can be related using the well-known formula for integration by parts \cite{CAB}, inverse functions \cite{KEY} or even Rolle's theorem \cite{CAR}. 

We wondered if there was an even deeper relation between the methods above and we have found that this is indeed the case. In particular we present here a method to compute the volume of a solid of revolution as a double integral in a very simple way. Then, we see that the classical methods (disks and shells) are recovered if this double integral is computed by each of the two possible applications of Fubini's theorem. As a further application we also show how Pappus' theorem is obtained from our formula.

\section{The volume as a double integral}
\label{doble}
\begin{figure}[ht]
\centering
\vspace{-0.35cm} \hspace{-1.0cm}
\includegraphics{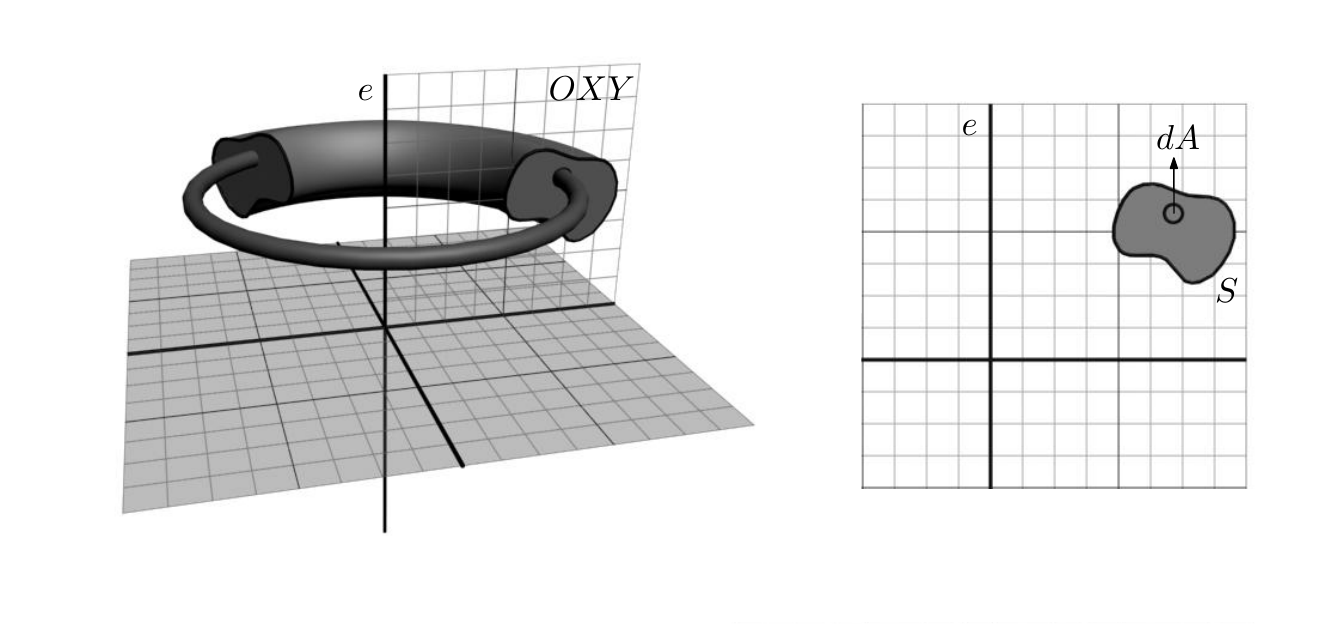}
\vspace{-1.5cm}
\caption{Solid of revolution together with the rotating section.}
\label{fig_ruben}
\end{figure}

Let $S$ be a closed region of the plane $OXY$ and let $e$ be any straight line in the same plane such that $e$ is exterior to $S$. For every point $P = (x,y)\in S$, put $d_e(x,y)=d(P,e)$ the distance from $P$ to $e$. Let us denote by $V(S,e)$ the volume of the solid obtained by rotating the region $S$ around the line $e$, see Figure \ref{fig_ruben}.

We claim that:
\begin{equation}\label{formula}
V(S,e)=\iint_S \, 2\pi d_e(x,y)\ dA.
\end{equation}

It is not the point here to give a complete and rigorous proof of this claim. The underlying idea is in fact very simple. For every point $P(x,y)\in S$, consider a tiny circle with center in $P$ and with area $dA$ (see the right hand side of Figure \ref{fig_ruben}). When this circle rotates around the axis $e$ it generates a torus of volume $2\pi d_e(x,y) dA$ and then it is enough to sum up all these volumes; i.e., to integrate over $S$.

Observe that $d_e(x,y)$ is always a polynomial of degree 1 in $x$ and $y$, namely, if the axis $e$ has equation $ax+by+c=0$, then $d_e(x,y)=\frac{|ax+by+c|}{\sqrt{a^2+b^2}}$.

\section{Deducing the ``standard'' methods}
\label{section_deducing}

In this section it is shown how the classical disk and shell methods can be obtained from the double integral formula presented above. 

\begin{figure}[ht]
\centering
\includegraphics{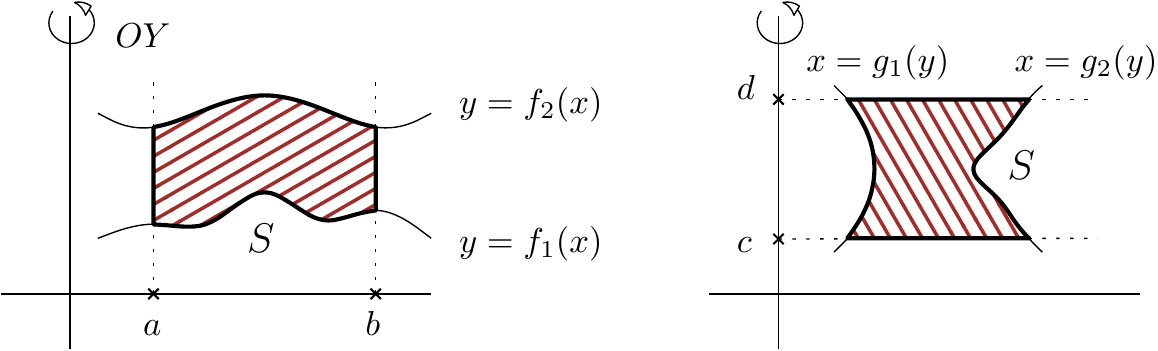}
\caption{Normal domains w.r.t. the $x$- and the $y$-axis, resp.}
\label{fig_deducing}
\end{figure}

Assume, without loss of generality, that the axis of revolution is $OY$. Let us first consider $S$ a normal domain with respect to the $x$-axis, i.e., the region $S$ is bounded by continuous functions $y=f_1(x)$ and $y=f_2(x)$ between $x=a$ and $x=b$ as in the left-hand side of Figure \ref{fig_deducing}. Fubini's theorem states that the double integral $I = \displaystyle{\iint_S \, 2\pi x \ dA}\,$ can be computed by means of simple integrals as:
$$
I = \int_{a}^{b} \left(\int_{f_1(x)}^{f_2(x)}2\pi x\ dy\right) \, dx = 
\int_a^b 2\pi \Big(f_2(x)-f_1(x) \Big) \, dx,
$$
which is precisely the well-known formula obtained in the method of shells.

Assume now that $S$ is a normal domain with respect to the $y$-axis. In this case $S$ is bounded by continuous functions $x=g_1(y)$ and $x=g_2(y)$ between $y=c$ and $y=d$ as in the right-hand side of Figure~\ref{fig_deducing}. Again, by Fubini, the double integral $I$ above can be solved as:
$$
\iint_S \, 2\pi x \ dA=\int_{c}^{d} \left(\int_{g_1(y)}^{g_2(y)}2\pi x\ dx\right) \, dy =
\int_c^d \pi \, \Big(g_2(y)^2 - g_1(y)^2 \Big) \, dy,
$$
which is the formula corresponding to the disk method.

This shows that, although both methods might seem very different geometrically (slicing the solid versus decomposing it in cylindrical shells), they are essentially the same. This is not surprising and has a simple interpretation. Using the formula $\displaystyle{\iint_S \, 2\pi x \ dA}\,$ one is adding up the volumes of the tori obtained by rotating a circle with center $(x,y)$ and area $dA$ about the $y$-axis. To sum all these volumes one can proceed in two different ways:
\begin{enumerate}
\item Fix $y=y_0$ and consider the sum of the volumes of the tori corresponding to the points $(x,y_0)$, obtaining a horizontal disk. Then it is enough to sum the volume of these disks. This is exactly the description of the disk method.
\item Fix $x=x_0$ and consider the sum of the volumes of the tori corresponding to the points $(x_0,y)$, obtaining a vertical cylindrical shell. Then one has just to sum the volumes of this shells. This is exactly the description of the shell method. 
\end{enumerate}
Hence the two classical methods arise from ours just by carefully arranging the tori whose volumes have to be summed.

\section{Pappus theorem}
\label{papus}
Consider a region $S$ of the plane $OXY$ and the solid obtained by rotating it around the axis $e$. If $C$ denotes the centroid of $S$ and $\mathcal{A}$ is the surface area of $S$ (recall the notation from Section \ref{doble}), then the so-called Pappus' theorem  states  in its classical form \cite[Chapter 6]{SHE} that the volume of this solid is given by
$$V(S,e)=2\pi d(C,e)\mathcal{A}.$$

Let us see how to obtain Pappus' theorem from Equation~\eqref{formula}. Recall that if $e$ has equation $ax+by+c=0$, then $d_e(x_0,y_0)=\frac{|ax_0+by_0+c|}{\sqrt{a^2+b^2}}$. It is also well-known that if $C=(x_C,y_C)$, then:
$$x_C=\frac{\iint_S x\ dA}{\mathcal{A}},\quad y_C=\frac{\iint_S y\ dA}{\mathcal{A}}.$$
With all this ingredients, and since $S$ can be assumed to lie on the semiplane determined by $ax+by+c>0$:
\begin{align*}
V(S,e)&=\iint_S 2\pi d_e(x,y)\ dA=2\pi\iint_S \frac{ax+by+c}{\sqrt{a^2+b^2}}\ dA=\\&=2\pi\mathcal{A}\frac{ax_C+by_C+c}{\sqrt{a^2+b^2}}=2\pi d(C,e)\mathcal{A},
\end{align*}
where linearity of the integral has played a key role.

\begin{figure}[ht]
\centering
\includegraphics{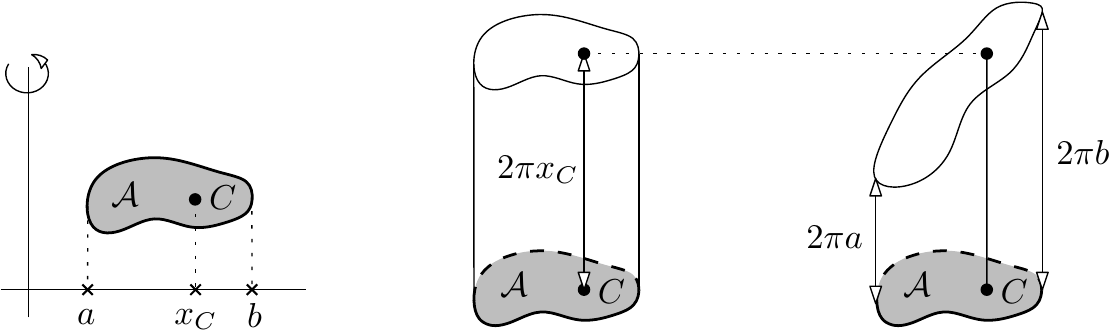}
\caption{Pappus' theorem \& our method.}
\label{fig_pappus}
\end{figure}

Pappus' theorem admits an interesting geometrical interpretation. When the region $S$ is rotated around $e$, a torus with section $S$ is obtained. Then Pappus' theorem implies that the volume of this torus is the same as that of a cylinder with base $S$ and height $2\pi d(C,e)$. In our case  we have a similar interpretation, but the volume of the original torus is now the same as that of a truncated cylinder. Figure~\ref{fig_pappus} illustrates both situations.

\section{An example}
\label{ejem}
As shown in Section \ref{section_deducing}, the formula for computing the volume of a solid of revolution using a double integral gives rise to the disk or shell method after applying Fubini's theorem. This fact could lead us to think that this approach does not provide in principle any new advantage with respect to the classical ones.


However, in order to solve a double integral, there are many more techniques at our disposal. The possibility to use specific double integration techniques, as for instance change of variables, is what gives merit to our method. An illustrative example based on changing to polar coordinates is shown.

\begin{figure}[ht]
\centering
\includegraphics{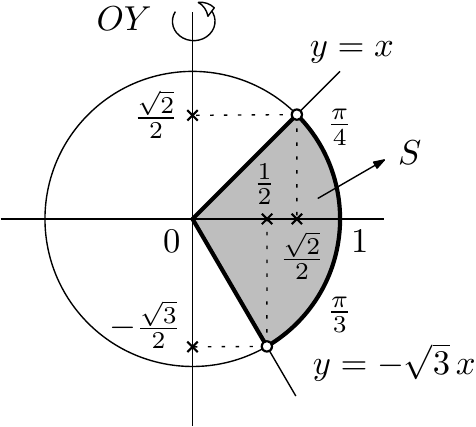}
\caption{Area of the rotating region about the $y$-axis.}
\label{fig_example}
\end{figure}

Let $S$ be the plane region bounded by the unit circle $x^2+y^2=1$ and the straight lines $y=x$ and $y=-\sqrt{3}\,x$\, as in Figure~\ref{fig_example}. Let us compute, using integrals, the volume of the solid obtained by rotating the region $S$ about the axis $OY$.

\begin{itemize}
\item Using the disk method, the required volume is:
\begin{equation}\label{f1}
V = \int_{-\frac{\sqrt{3}}{2}}^{0} \ \pi \, \Big( 1-y^2-\frac{y^2}{3} \Big) \, dy + \int_0^{\frac{\sqrt{2}}{2}} \pi \, \Big( 1-y^2-y^2 \Big) \, dy.
\end{equation}
\item Using the method of shells, the volume is:
\begin{equation}\label{f2}
\begin{split}
V = & \int_{0}^{\frac{1}{2}} 2 \pi x \, \Big( x+\sqrt{3}x \Big) \, dx + \int_{\frac{1}{2}}^{\frac{\sqrt{2}}{2}} 2 \pi x \, \Big( x-\sqrt{1-x^2} \ \Big) \, dx \\
& + \int_{\frac{\sqrt{2}}{2}}^{1} \ 2\pi \, \Big( 2\sqrt{1-x^2} \ \Big) \, dx.
\end{split}
\end{equation}
\item Finally, using our method, the volume is:
\begin{equation}
V=\iint_S \, 2 \pi x \ dA.
\end{equation}
\end{itemize}

Observe that the last expression using double integrals is much simpler and more compact. Applying Fubini's theorem one recovers the two other expressions and there would be no real advantage in our method. Instead, let us change to polar coordinates.

Then,
$$
V = 2\pi\int_{0}^{1}\rho^2 \, d\rho \ \int_{-\frac{\pi}{3}}^{\frac{\pi}{4}}\cos\theta\, d\theta=2\pi\left[\frac{\rho^3}{3}\right]_0^{1} \, \Big[\sin\theta\Big]_{-\frac{\pi}{3}}^{\frac{\pi}{4}}=\frac{\pi(1+\sqrt{2})}{3},$$
and thus the required volume can be found in a simple way compared to the expressions~\eqref{f1} and~\eqref{f2}. 

\section{Concluding remarks}
\label{conc}
We have presented a way to compute the volume of a solid of revolution as a double integral. This method seems, as far a we know, to be absent in the literature. In addition to its intrinsic mathematical interest, we think that this method to compute the volume of a solid of revolution might be of educational interest. Specifically, some of its advantages are:
\begin{enumerate}
\item It avoids considerations about the shape of the solid.
\item It gives an easy way to describe the volume of the solid when the axis of revolution is not horizontal neither vertical.
\item It introduces the use of double integration techniques (polar coordinates, for instance) that, in some cases (see the example in Section 5 above), leads to shorter calculations.
\item It extends the disk and shell methods and also Pappus' theorem in the sense that all of them can easily be recovered from Equation~\eqref{formula}.
\end{enumerate}

\noindent {\bf Acknowledgments.} We would like to thank Rubén Vigara for his help with Figure 1.

\end{document}